\documentclass[leqno,12pt]{article}

\usepackage[dvipdfmx]{graphicx}
\usepackage{amsmath}
\usepackage{txfonts}
\usepackage{fancyhdr}
\usepackage{latexsym}
\title{Non-orientable genus of a knot in~punctured~$\mathbb{C}P ^2$}
\author{KOUKI SATO AND MOTOO TANGE}
\date{}

\newtheorem{dfn}{Definition}

\newtheorem{problem}{Problem}
\newtheorem{thm}{Theorem}[section]
\newtheorem{prop}[thm]{Proposition}
\newtheorem{lem}{Lemma}[section]

\begin{document}
\maketitle
\thispagestyle{empty}

\begin{abstract}
For any knot $K$ which bounds non-orientable and null-homologous
surfaces $F$ in punctured $n\mathbb{C}P^2$,
we construct a lower bound of the first Betti number of $F$ 
which consists of the signature of $K$ and the Heegaard Floer $d$-invariant 
of the integer homology sphere obtained by $1$-surgery along $K$.
By using this lower bound, we prove that for any integer $k$, a certain knot
cannot bound any surface which satisfies the above conditions and whose first Betti number is less than $k$.
\end{abstract}

\section{Introduction}
\paragraph{}
Throughout this paper, we work in the smooth category, all 4-manifolds
are orientable, oriented and simply-connected, and all surfaces are compact. 
If $M$ is a closed 4-manifold, punc$M$ denotes $M$ with an open 4-ball
deleted.

The non-orientable 4-genus $\gamma_4 (K)$ of a knot $K$ 
is the smallest first Betti number of any non-orientable surface
in $B^4$ with boundary $K$. 
It has been investigated in \cite{batson}, \cite{gilmer-livingston},
\cite{murakami-yasuhara}, \cite{viro}, and \cite{yasuhara}.
In this paper, we extend the definition of $\gamma_4 (K)$ to any 4-manifold with boundary
$S^3$.
\begin{dfn}
Let $M$ be a closed 4-manifold and
$K \subset \partial (\text{punc}M)\  (\cong S^3)$ a knot.
The non-orientable $M$-genus  $\gamma_M (K)$ of $K$
is the smallest first Betti number of any non-orientable
surface $F \subset \text{punc}M$ with boundary $K$.

Moreover, we define $\gamma^0_{M}(K)$ 
to be 
the smallest first Betti number of any non-orientable
surface $F \subset \text{punc}M$ with boundary $K$
which  represents zero in 
$H_2 (\text{punc}M, \partial(\text{punc}M); \mathbb{Z}_2 )$.
\end{dfn}
We note that $\gamma_{4}(K) = \gamma_{S^4}(K) = \gamma^0_{S^4}(K)$ and
$\gamma_{M}(unknot) = \gamma^0_{M}(unknot)=1$ for any closed 4-manifold $M$.
In this paper, we  consider the following problem.
\begin{problem}
Can $\gamma_{M}$ and $\gamma^0_{M}$ be taken arbitrarily large?
\end{problem}
The answer of this problem depends on the choice of $M$.
For example, Suzuki \cite{suzuki} proved that 
any knot bounds a disk in punc($S^2 \times S^2$) and 
 punc($\mathbb{C}P^2 \# \overline{\mathbb{C}P^2}$).
It follows that $\gamma_{S^2 \times S^2}(K) = \gamma_{\mathbb{C}P^2 \# \overline{\mathbb{C}P^2}}(K)
= 1$ for any knot $K$.
For a long time, it had been unknown whether Problem 1 even in the case $\gamma_{S^4}$ is true or not.
Most recently Batson in \cite{batson} gave the affirmative answer for the problem by using the following inequality:
\begin{thm}[Batson,\cite{batson}]
\label{thm.b}
Let $K \subset S^3$ be a knot. Then
\[
\gamma_{S^4}(K) \geq \frac{-\sigma (K)}{2} + d(S^3_1(K)),
\]
where $\sigma$ denotes the signature of $K$ and $d(S^3_1(K))$ the Heegaard-Floer
$d$-invariant of the integer homology sphere obtained by $1$-surgery along $K$.
\end{thm}

The $d$-invariant is defined by Ozsv\'{a}th and Szab\'{o} in \cite{ozsvath-szabo}.
We extend Theorem \ref{thm.b} to the case of $\gamma^0_{n\mathbb{C}P2}$. 

\begin{thm}
\label{thm2}
Let $K \subset S^3$ be a knot. Then
\[
\gamma^0_{n\mathbb{C}P^2}(K) \geq \frac{-\sigma (K)}{2} + d(S^3_1(K)) -n.
\]
\end{thm}

By applying Theorem \ref{thm2}, we give the answer of Problem 1 in the case of 
$\gamma^0_{n\mathbb{C}P^2}$.

\begin{thm}
\label{thm3}
For every $k$,
there exists a knot $K$ such that
$
\gamma^0 _{n\mathbb{C} P^2  } (K) = k.
$
\end{thm}

In fact, we show that 
$\gamma^0 _{n\mathbb{C} P^2  } (\overset{n+k}{\#}9_{42}) = k$, where
$\overset{n+k}{\#}9_{42}$ denotes the $n+k$ times connected sum of  
the knot $9_{42}$ in Rolfsen's table \cite{rolfsen}. 
In order to prove Theorem \ref{thm2}, we first prove the following proposition.

\begin{prop}
\label{thm1}
Let $K \subset \partial (\text{punc}(n\overline{\mathbb{C}P^2}))$ be a knot
and $F \subset punc(n\overline{\mathbb{C}P^2})$ a non-orientable
surface with boundary $K$. Then
\begin{equation*}
\beta_1(F) \geq \frac{e(F)}{2} - 2d(S^3_{-1}(K)),
\end{equation*}
where $e(F)$ is the normal Euler number of $F$,
and $\beta_{i}$ denotes the $i$-th Betti number.
\end{prop}

This proposition is an extension of  Theorem $4$ in \cite{batson}
and independent of the homology class of $F$.
However, the normal Euler number $e(F)$ depends on
$F$, hence this lower bound is not
enough to determine the genus.
To delete $e(F)$, we use the following theorem.
\begin{thm}[Yasuhara, \cite{yasuhara} ]
\label{thm.y}
Let $M$ be a closed 4-manifold and
$K \subset \partial (\text{punc}M)$ a knot.
If $K$ bounds a non-orientable surface $F$ in punc$M$
that represents zero in $H_2 (\text{punc}M, \partial (\text{punc}M); \mathbb{Z}_2 ) $, then
\[
\left| \sigma (K) + \sigma (M) -  \frac{e(F)}{2} \right| \leq \beta_2 ( M ) + \beta_1 ( F ),
\]
where $\sigma (M)$ is the signature of $M$.
\end{thm}

By applying this theorem to $n\mathbb{C}P^2$ and deleting the term of $e(F)/2$,
we obtain the inequality of Theorem \ref{thm2}.

\paragraph{Remark.}
In our arguments in Section~\ref{ProofofLem1.1} and Section~\ref{proofofthm1.2}, it is clear that the inequalities in
Theorem~\ref{thm2} and Proposition~\ref{thm1} hold for any 4-manifold $P$ with intersection form $n\langle 1\rangle$ instead of $n{\mathbb C}P^2$,
because these homological data of the ambient space determines the inequalities.
For simplicity we deal with the case of $n{\mathbb C}P^2$.
\paragraph{Acknowledgements.}
The first author would like to thank Kokoro Tanaka and Akira Yasuhara 
for their useful comments and encouragement.

\section{Proof of Proposition \ref{thm1}}

In order to prove Proposition \ref{thm1}, we first prove the following lemma.

\begin{lem}
\label{prop2.1}
Under the hypothesis of Lemma \ref{thm1}, moreover if $\beta_1(F)$ is odd, then
there exists an orientable surface
$
F' \subset punc(n\overline{\mathbb{C} P^2} \  \sharp \  S^2 \times S^2 ) 
$
which is still bounded by $K$, and has $\beta_1 (F') = \beta_1 (F) - 1 $
and $e(F') = e(F) + 2 $.
\end{lem}  

\begin{figure}
\begin{center}
\includegraphics[clip]{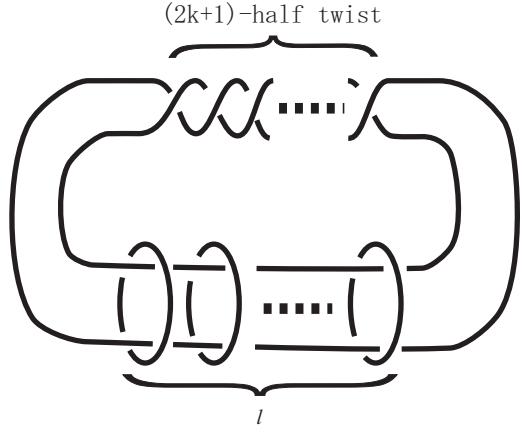}\\
\caption{the link $L$ (for $k \geq 0$)}
\label{link L}
\end{center}
\end{figure}

\paragraph{Proof.}
Since $\beta_1(F)$ is odd, there is a curve $C$ in $F$ such that 
the regular neighborhood of $C$ in F is a M\"{o}bius band and
$F \setminus C$ is orientable.
By the simply-connectedness of punc$(n\overline{\mathbb{C} P^2 })$,
 $C$ is null-homotopic. We note that in these dimensions (i.e.,
for 1-manifolds in 4-manifolds) every homotopy may be replaced 
with an isotopy. It follows that $C$ bounds an embedded 2-disk $D$ in 
punc$(n\overline{\mathbb{C}P^2})$.
Without loss of generality we can assume that $D$ is transverse to $F$.
Then $F \cap D$ is the disjoint union of the curve $C$ and some 
transversal intersections $\{p_i\}$ $(i = 1,2,...,l)$. 
Let  $V(D)$ be a small regular neighborhood of $D$ in punc$(n\overline{\mathbb{C}P^2})$.
$V(D)$ is diffeomorphic to $D \times D^2$, and $F\cap V(D)$ consists of one M\"{o}bius band 
and $l$ $2$-disks $p_i \times D^2 $. 
If we draw $\partial V(D) \cong S^3 $ with its standard decomposition into
solid tori 
$
\partial V(D) \cong 
\partial D \times D^2 \cup_{\partial D \times S^1} D \times S^1 
$,
we see $\partial(F \cap V(D)) \subset \partial V(D)$ as the link $L$
in Figure~{\ref{link L}} consisting of a
$(2, 2k + 1)$-cable of the core for the first factor and
$l$ parallel copies of the core for the second. 
It was proved by \cite{batson} and \cite{yasuhara} that $L$ bounds $l+1$ disjoint embedded disks $E$ in punc$(S^2 \times S^2)$
such that $e(E) = e(F \cap V(D)) + 2$.
We remark that $[E,\partial E] \in H_2(\text{punc}(S^2 \times S^2)
;\mathbb{Z})$ is $2\alpha + b\beta$ $(b \in \mathbb{Z})$,
where $\alpha$ and $\beta$ are standard generators of 
$
H_2 (\text{punc}(S^2 \times S^2)
, \partial(\text{punc}(S^2 \times S^2 ))
;\mathbb{Z})
$
such that
$\alpha \cdot \alpha = \beta \cdot \beta = 0$, and $\alpha \cdot \beta = 1$.
Let $F'' = F \setminus (F \cap V(D))$.
Excising $V(D)$ from punc($n\overline{\mathbb{C} P^2}$) and 
capping off the pair ($\partial V(D)$, $L$) with a pair 
(punc($S^2 \times S^2$), $E$),  we obtain a new orientable surface
$F' = F'' \cup E$ in 
punc($n\overline{\mathbb{C} P^2} \  \sharp \  S^2 \times S^2 $)
with boundary $K$.
It is easy to check that $\beta_1 (F') = \beta_1 (F) - 1$,
and by the additivity of the normal Euler number, $e(F') = e(F) + 2$. 
This completes the proof.
\hspace{\fill} $\Box$ \\

We note that  the homology class
$
[F', \partial F']  \in 
H_2 (\text{punc}(n\overline{\mathbb{C} P^2} \  \sharp \  S^2 \times S^2 )
, \partial(\text{punc}(n\overline{\mathbb{C} P^2} \  \sharp \  S^2 \times S^2 ))
;\mathbb{Z})
$
is 
\[
\sum^{j}_{i=1}2a_i \overline{\gamma_i} + 
\sum^{n}_{i= j+1}(2a_{i}+1) \overline{\gamma_i} + 2\alpha + b\beta
\ \  (a_i, j \in \mathbb{Z}, 0 \leq j \leq n ),
\]
where $\overline{\gamma_i}$, are standard generators of 
$
H_2 (\text{punc}(n\overline{\mathbb{C} P^2})
, \partial(\text{punc}(n\overline{\mathbb{C} P^2}))
;\mathbb{Z})
$
such that $\overline{\gamma_i} \cdot \overline{\gamma_j} = -\delta_{ij}$ (Kronecker's delta).
Since $F'$ is orientable, 
$
e(F') = [F',\partial F'] \cdot [F',\partial F'] = 
-\sum^{j}_{i=1}4a^2_i -\sum^{n}_{i= j+1}(2a_{i}+1)^2 + 4b.
$

We next prove the following lemma. It is a generalization of a discussion in 
Section  3 \cite{batson}.

\begin{lem}
\label{new lem}
Let $M$ be an integer homology 3-sphere,
$X$ a simply-connected 4-manifold such that $\partial X = M$ and $\beta_2^+(X)=1$
and
$\Sigma$ an orientable closed surface in $X$ with genus $g$
and self-intersection $m$.
Then for any Spin$^c$ structure $\mathbf{s}$ of $X$ which satisfies
$\left<c_1( \mathbf{s}),[\Sigma]\right> = m-2g > 0$,
the following inequality holds:
\[
c_1( \mathbf{s})^2 + \beta^-_2(X) \leq 1 + 4d(M).
\]
\end{lem}

\paragraph{Proof of Lemma \ref{new lem}}
Suppose that $ X' = X \setminus \nu (\Sigma ) $, then $X'$
is a negative semi-definite 4-manifold 
with disconnected boundaries $Y_{g,-m} \amalg M$,
where $Y_{g,-m} $ denotes the Euler
number $-m$ circle bundle over $\Sigma$.
We apply the following theorem to the pair $(X', Y_{g,-m} \amalg M)$.

\begin{thm}[Ozsv\'{a}th and Szab\'{o}, \cite{ozsvath-szabo}]
\label{thm.o-s}
Let $Y$ be a closed oriented 3-manifold (not necessarily connected)
with standard $HF^{\infty} $, endowed with a torsion Spin$^c$ structure 
$\mathbf{t}$.
If $X$ is a negative semi-definite four-manifold bounding $Y$ such that
the restriction map $H^1 (X;\mathbb{Z}) \rightarrow H^1 (Y;\mathbb{Z})$ is
trivial, and $\mathbf{s}$ is a Spin$^c$ structure on $X$ restricting to $\mathbf{t}$
on $Y$, then
\begin{equation*}
c_1 (\mathbf{s} )^2 + \beta^- _2 (X) \leq 4d_b (Y,\mathbf{t} ) + 2\beta_1 (Y).
\end{equation*}
\end{thm}
It  was proved by \cite{batson} and \cite{ozsvath-szabo} that 
$Y_{g,-m} \amalg S^3_{-1} (K)$
has standard $HF^{\infty}$.
Since we can verify that $H^1(X')=0$
in the same way as Section 3 in \cite{batson},
it follows that the pair $(X', Y_{g,-m} \amalg M)$ satisfies all conditions of
Theorem \ref{thm.o-s}.

Applying Theorem \ref{thm.o-s} to the Spin$^c$ structure $\mathbf{s}_t|_{X'}$ on $X'$,
we have
\begin{equation}
c_1 (\mathbf{s} |_{X'} )^2 + \beta^- _2 (X') \leq 4d_b (Y_{g,-m}, \mathbf{s} |_{Y_{g,-m}})
+ 4d(M) + 2\beta_1 (Y_{g,-m}) + 2\beta_1(M).
\end{equation}
Let us compute each term in this inequality (1).
In order to compute $c_1(\mathbf{s}_t|_{X'})^2$, we decompose the intersection
form of $X$ in terms of the $\mathbb{Q}$-valued intersection forms on
$\nu (\Sigma )$ and $X'$ ; if $c \in H^2 (X)$, then
\begin{equation*}
Q_{X} (c) = Q_{\nu (\Sigma )} (c|_{\nu (\Sigma)})+ Q_{X'} (c|_{X'}).
\end{equation*}
This gives 
$c_1(\mathbf{s})^2 = c_1(\mathbf{s}|_{\nu(\Sigma)})^2 + c_1(\mathbf{s}|_{X'})^2$.
Hence we have
\begin{equation*}
c_1 (\mathbf{s} |_{X'} )^2 = c_1 (\mathbf{s})^2 -
c_1(\mathbf{s}|_{\nu(\Sigma)})^2 =
c_1 (\mathbf{s})^2 - \frac{(m - 2g)^2}{m}.
\end{equation*}
For the above Spin$^c$ structure $\mathbf{s}|_{Y_{g,-m}}$,
the $d$-invariant of $Y_{g,-m}$ is computed in section 9
of \cite{ozsvath-szabo}.
If $\langle c_1(\mathbf{s} ), [\Sigma] \rangle =m -2g > 0$, then
\begin{equation*}
d_b (Y_{g,-m}, \mathbf{s} |_{Y_{g,-m}}) = \frac{1}{4} - \frac{g^2}{m} - \frac{m}{4}.
\end{equation*}
After substituting all the values computed above, (1) reduces to
\begin{equation}
c_1 (\mathbf{s})^2 - \frac{(m - 2g)^2}{m}
 + \beta^- _2 (X') 
\leq 4 \left(\frac{1}{4} - \frac{g^2}{m} - \frac{m}{4} \right)
+ 4d(S^3_{-1} (K)) + 2(2g).
\end{equation}
Since $\beta^- _2 (X')=\beta^-_2(X)$, (2) gives the inequality
\[
c_1( \mathbf{s})^2 + \beta^-_2(X) \leq 1 + 4d(M).
\]
\hspace{\fill} $\Box$ \\

\paragraph{Proof of Proposition \ref{thm1}.}
\label{ProofofLem1.1}
Note that for any knot $K$, $d(S^3_{-1}(K)) \geq 0$.
Hence when $e(F) \leq \beta_1(F) $, it is clear that
this proposition holds.
Therefore we assume that $e(F) > \beta_1(F) $.

We first give the proof for the case where $\beta_1(F)$ is odd. 
In Lemma \ref{prop2.1}, we constructed an orientable surface 
$
F' \subset (\text{punc}(\overline{\mathbb{C} P^2} \ \sharp \  S^2 \times S^2 ))
$
with boundary $K \subset S^3$.
Attaching a ($-1$)-framed 2-handle along $K$,  we have a 4-manifold $\overline{W}$
with boundary $S^3_{-1} (K) $ and the intersection form
\begin{equation*}
Q_{\overline{W}} =
\left(
\begin{array}{c|ccc|cc}
-1&0&\ldots&0&0&0\\ \hline
0&-1&\ &O&0&0\\
\vdots&\ &\ddots&\ &\vdots&\vdots\\
0&O&\ &-1&0&0\\ \hline
0&0&\ldots&0&0&1\\
0&0&\ldots&0&1&0
\end{array}
\right).
\end{equation*}
We may cap off $F'$ with the core of the 2-handle to form a closed
surface $\Sigma$ with genus $g = (b_1 (F) - 1)/2$, homology class
$
\overline{\gamma_0} +
\sum^{j}_{i=1}2a_i \overline{\gamma_i} + 
\sum^{n}_{i= j+1}(2a_{i}+1) \overline{\gamma_i} + 2\alpha + b\beta
$, and the self-intersection number
\begin{equation*}
m = -1 -\sum^{j}_{i=1}4a^2_i -\sum^{n}_{i= j+1}(2a_{i}+1)^2 + 4b= e(F) + 1 > 0.
\end{equation*}

We next choose a Spin$^c$ structure on $W$.
Since $H^2 (\overline{W}) \cong \mathbb{Z}^{m+3}$ has no 2-torsion,
a Spin$^c$ structure on $\overline{W}$ is determined by its first Chern class.
Fix a Spin$^c$ structure $\mathbf{s}_t$ on $\overline{W}$ satisfying
\begin{equation*}
PD(c_1 (\mathbf{s}_t ) ) =
\varepsilon \overline{\gamma_0} +
\sum^{n}_{i= 1}(2a_{i}+1) \overline{\gamma_i} + 2\alpha + 2x\beta,
\end{equation*}
where
\begin{equation*}
x = \frac{\sum^{j}_{i= 1}2a_{i} + 2(b - g ) - 1 + \varepsilon }{4}
\end{equation*}
and $\varepsilon \in \{ 1, -1 \}$  is chosen so as to make $x$ an integer.
Since the given vector is characteristic for $Q_{\overline{W} }$, it
corresponds to a Spin$^c$ structure. Furthermore, 
$\langle c_1(\mathbf{s}_t ), [\Sigma] \rangle = m - 2g = e(F) - \beta_1(F)+2>0$.
Applying Lemma \ref{new lem} to the pair $(\overline{W}, S^3_{-1}(K))$,
we have
\begin{equation}
c_1( \mathbf{s}_t)^2 + \beta^-_2(\overline{W}) \leq 1 + 4d(S^3_{-1}(K)).
\end{equation}
Since 
$
c_1 (\mathbf{s}_t)^2 =
-1 - \sum^{n}_{i=1} (2a_i+1)^2 + 8x 
= e(F) -j -1 +2\varepsilon -4g,
$
(3) gives
\begin{equation}
\left(e(F) -j -1 +2\varepsilon - 4g \right) + (n+2) \leq 1+ 4d(S^3_{-1} (K)).
\end{equation}
Since $-1 \leq \varepsilon$, $j \leq n$, and $2g = \beta_1(F) - 1$, 
(4) reduces to the following inequality
\begin{equation}
\frac{e(F)}{2} - 2d \big( S ^3 _{-1} (K) \big)
\leq \beta_1 (F).
\end{equation}

Finally, we consider the case where $\beta_1(F)$ is even.
Taking the connected sum 
$F \subset \text{punc}(n\overline{\mathbb{C}P^2})$
with the standard embedding of $\mathbb{R}P^2 \subset S^4$ whose 
normal Euler number
is $+2$, and we have a  non-orientable surface 
$\hat{F} \subset \text{punc}(n\overline{\mathbb{C}P^2})$
with boundary $K$
such that $\beta_1 (\hat{F}) = \beta_1 (F) + 1 $
and $e(\hat{F}) = e(F) + 2 $. 
Since $\beta_1(\hat{F})$ is odd, 
$\hat{F}$ satisfies the inequality (5).
This implies that the inequality (5) holds in this case.

This completes the proof of Proposition \ref{thm1}. 
\hspace{\fill} $\Box$ \\

\section{Proof of Theorem \ref{thm2}}
\label{proofofthm1.2}
By reversing the orientation of $M$, we obtain the following proposition. 
We use this proposition to prove Theorem \ref{thm2}.
\begin{prop}
\label{prop2.1}
For any 4-manifold $M$ and any knot $K$, the following equality holds;
\[
\gamma^0_{M}(K)=\gamma^0_{-M}(\overline{K}),
\]
where $\overline{K}$ denotes the mirror image of $K$.
\end{prop}

Let us prove Theorem \ref{thm2}.

\paragraph{Proof of Theorem \ref{thm2}.}

Suppose that $F  \subset \text{punc}(n\overline{\mathbb{C}P^2})$
is a non-orientable surface 
with boundary $K$
which represents zero in 
$
H_2 (\text{punc}(n\overline{\mathbb{C}P^2}),
\partial(\text{punc}(n\overline{\mathbb{C}P^2})), \mathbb{Z}_2 )
$.
It follows from Theorem \ref{thm.y} that
\[
\left| \sigma (K) + (-n) -  \frac{e(F)}{2} \right| \leq n + \beta_1 ( F ).
\]
Hence we have
\[
\beta_1 ( F ) \geq \sigma (K) -  \frac{e(F)}{2} - 2n.
\]
Combining this inequality with Proposition \ref{thm1}, we have
\[
\gamma^0_{n\overline{\mathbb{C}P^2}}(K) \geq \frac{\sigma (K)}{2} - d(S^3_{-1} (K)) - n.
\]
By using this inequality and Proposition \ref{prop2.1}, it follows that for any knot 
$K \subset \partial ( \text{punc} (n\mathbb{C}P^2))$,
\[
\gamma^0_{n\mathbb{C}P^2}(K) =
\gamma^0_{n\overline{\mathbb{C}P^2}}(\overline{K}) 
\geq \frac{\sigma (\overline{K})}{2} - d(S^3_{-1} (\overline{K})) - n =
 \frac{-\sigma (K)}{2} + d(S^3_{1} (K)) - n.
\]
This proves Theorem \ref{thm2}.
\hspace{\fill} $\Box$ \\

\section{Proof of Theorem \ref{thm3}}

To prove Theorem~\ref{thm3}, it is necessary to show $d(S^3_1(\overset{m}{\#}K))=0$ for some knot $K$.
Since $d(S^3_1(\cdot))$ is a knot concordance invariant by the result of \cite{peters}, but not a homomorphism, we must prove $d(S^3_1(\overset{m}{\#}K))=0$ for each $m$ and
the knot $K$. 
Throughout this paper, the coefficient $\mathbb{F}$ of any Heegaard Floer homology is the field with the order $2$.
The coordinate $(i,j)$, as is used below, is the same as that in \cite{ozsvath-szabo2}.
We denote the whole differential in the knot Floer chain complex $CFK^\infty(K)$ by $\partial^\infty$ and
denote the differential restricted to vertical (or horizontal) lines by $\partial^{vert}$ (or $\partial^{hor}$ respectively).
For the other differentials, we use the same notation $\partial^\infty$.
\begin{prop}
\label{tangeprop}
For any positive integer $m$, $d(S^3_1(\overset{m}{\#}9_{42}))=0$.
\end{prop}
\paragraph{Proof of Proposition \ref{tangeprop}.}
Due to \cite{peters}, the correction term $d(S^3_1(K))$ coincides with $\tilde{d}(S^3_p(K),[0])$, where $p$ is a sufficient large
integer.
The correction term $\tilde{d}$ is the unshifted correction term for $CFK^\infty(K)\{\max(i,j)\ge 0\}$,
namely $\tilde{d}(S^3_p(K),[0])=d(S^3_p(K),[0])-\frac{p-1}{4}$.

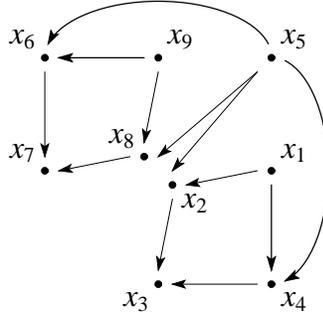
\begin{figure}[htbp]
\begin{center}
\unitlength 0.1in
\begin{picture}( 20.6300, 17.9900)(  9.3700,-28.0000)
%
\special{pn 20}%
\special{sh 1}%
\special{ar 1432 1298 10 10 0  6.28318530717959E+0000}%
\special{sh 1}%
\special{ar 2026 1298 10 10 0  6.28318530717959E+0000}%
\special{sh 1}%
\special{ar 2618 1298 10 10 0  6.28318530717959E+0000}%
\special{sh 1}%
\special{ar 2618 1890 10 10 0  6.28318530717959E+0000}%
\special{sh 1}%
\special{ar 2618 2484 10 10 0  6.28318530717959E+0000}%
%
\special{pn 8}%
\special{pa 1952 1298}%
\special{pa 1506 1298}%
\special{fp}%
\special{sh 1}%
\special{pa 1506 1298}%
\special{pa 1574 1318}%
\special{pa 1560 1298}%
\special{pa 1574 1278}%
\special{pa 1506 1298}%
\special{fp}%
%
\special{pn 8}%
\special{pa 2618 1964}%
\special{pa 2618 2410}%
\special{fp}%
\special{sh 1}%
\special{pa 2618 2410}%
\special{pa 2638 2342}%
\special{pa 2618 2356}%
\special{pa 2598 2342}%
\special{pa 2618 2410}%
\special{fp}%
%
\special{pn 8}%
\special{ar 2026 1298 594 296  3.3893469 6.0354311}%
%
\special{pn 8}%
\special{pa 1458 1212}%
\special{pa 1450 1224}%
\special{fp}%
\special{sh 1}%
\special{pa 1450 1224}%
\special{pa 1500 1176}%
\special{pa 1476 1178}%
\special{pa 1464 1156}%
\special{pa 1450 1224}%
\special{fp}%
%
\special{pn 8}%
\special{ar 2618 1890 298 594  4.9565753 6.2831853}%
\special{ar 2618 1890 298 594  0.0000000 1.3266101}%
%
\special{pn 8}%
\special{pa 2704 2458}%
\special{pa 2692 2466}%
\special{fp}%
\special{sh 1}%
\special{pa 2692 2466}%
\special{pa 2760 2452}%
\special{pa 2738 2440}%
\special{pa 2740 2416}%
\special{pa 2692 2466}%
\special{fp}%
\put(26.6300,-12.5200){\makebox(0,0)[lb]{$x_5$}}%
\put(26.6300,-18.4500){\makebox(0,0)[lb]{$x_1$}}%
\put(20.6900,-12.5200){\makebox(0,0)[lb]{$x_9$}}%
%
\special{pn 20}%
\special{sh 1}%
\special{ar 1952 1816 10 10 0  6.28318530717959E+0000}%
\special{sh 1}%
\special{ar 2100 1964 10 10 0  6.28318530717959E+0000}%
%
\special{pn 20}%
\special{sh 1}%
\special{ar 1432 1890 10 10 0  6.28318530717959E+0000}%
\special{sh 1}%
\special{ar 2026 2484 10 10 0  6.28318530717959E+0000}%
%
\special{pn 4}%
\special{pa 2026 1372}%
\special{pa 1952 1742}%
\special{fp}%
\special{sh 1}%
\special{pa 1952 1742}%
\special{pa 1984 1680}%
\special{pa 1962 1690}%
\special{pa 1944 1672}%
\special{pa 1952 1742}%
\special{fp}%
%
\special{pn 4}%
\special{pa 2544 1890}%
\special{pa 2174 1964}%
\special{fp}%
\special{sh 1}%
\special{pa 2174 1964}%
\special{pa 2242 1972}%
\special{pa 2226 1954}%
\special{pa 2234 1932}%
\special{pa 2174 1964}%
\special{fp}%
%
\special{pn 4}%
\special{pa 2100 2038}%
\special{pa 2026 2410}%
\special{fp}%
\special{sh 1}%
\special{pa 2026 2410}%
\special{pa 2058 2348}%
\special{pa 2036 2358}%
\special{pa 2018 2340}%
\special{pa 2026 2410}%
\special{fp}%
%
\special{pn 4}%
\special{pa 1878 1816}%
\special{pa 1506 1890}%
\special{fp}%
\special{sh 1}%
\special{pa 1506 1890}%
\special{pa 1576 1898}%
\special{pa 1558 1880}%
\special{pa 1568 1858}%
\special{pa 1506 1890}%
\special{fp}%
%
\special{pn 4}%
\special{pa 2544 2484}%
\special{pa 2100 2484}%
\special{fp}%
\special{sh 1}%
\special{pa 2100 2484}%
\special{pa 2166 2504}%
\special{pa 2152 2484}%
\special{pa 2166 2464}%
\special{pa 2100 2484}%
\special{fp}%
\special{pa 1432 1372}%
\special{pa 1432 1816}%
\special{fp}%
\special{sh 1}%
\special{pa 1432 1816}%
\special{pa 1452 1750}%
\special{pa 1432 1764}%
\special{pa 1412 1750}%
\special{pa 1432 1816}%
\special{fp}%
%
\special{pn 4}%
\special{pa 2544 1372}%
\special{pa 2026 1816}%
\special{fp}%
\special{sh 1}%
\special{pa 2026 1816}%
\special{pa 2090 1788}%
\special{pa 2066 1782}%
\special{pa 2064 1758}%
\special{pa 2026 1816}%
\special{fp}%
%
\special{pn 4}%
\special{pa 2544 1372}%
\special{pa 2100 1890}%
\special{fp}%
\special{sh 1}%
\special{pa 2100 1890}%
\special{pa 2158 1852}%
\special{pa 2134 1850}%
\special{pa 2128 1826}%
\special{pa 2100 1890}%
\special{fp}%
\put(26.6300,-25.2700){\makebox(0,0)[lt]{$x_4$}}%
\put(13.8700,-12.5200){\makebox(0,0)[rb]{$x_6$}}%
\put(13.8700,-18.4500){\makebox(0,0)[rb]{$x_7$}}%
\put(19.0600,-17.7100){\makebox(0,0)[rb]{$x_8$}}%
\put(21.4400,-20.0800){\makebox(0,0)[lt]{$x_2$}}%
\put(19.8100,-25.2700){\makebox(0,0)[rt]{$x_3$}}%
%
\special{pn 8}%
\special{pa 3000 2800}%
\special{pa 1200 2800}%
\special{ip}%
\end{picture}%
\caption{The differential maps of the fundamental part $G_0$ of $CFK^\infty(9_{42})$ (Fig.14 in \cite{ozsvath-szabo2}).}
\label{complex942}
\end{center}
\end{figure}

For the generators of $CFK^\infty(9_{42})$, we use the same as those in Fig.14 in \cite{ozsvath-szabo2} (see Figure~\ref{complex942}).
Let denote $S_1=\{x_i\,|\,1\le i\le 9\}$.
Let $G_i$ be a graded module $\mathbb{F}\langle U^{-i}\cdot x\,|\,x\in S_1\rangle$, where $x_5$ has the Alexander grading $gr(x_5)=0$.
We call the chain complex $G_0$ a {\it fundamental part} in $CFK^\infty$.
$U$ is the action decreasing the grading by $2$.
The chain complex $CFK^\infty(9_{42})$ consists of an infinite direct sum $\oplus_{i\in\mathbb{Z}}G_i$
and the class $\alpha_i:=U^{-i}(x_1+x_5+x_9)$ is the homological generator of $CFK^\infty(9_{42})$.
By the result in \cite{ozsvath-szabo2}, for the large $p$-Dehn surgery we obtain the
graded isomorphism
$$HF^+_{\ast+\frac{p-1}{4}}(S^3_p(9_{42}),[0])\cong H_\ast(CFK^\infty(9_{42})\{\max(i,j)\ge 0\}).$$

The right hand side is easily computed as $H_\ast(CFK^\infty(9_{42})\{\max(i,j)\ge 0\})=\oplus_{i\ge 0}\mathbb{F}\cdot\alpha_i$ and the minimal grading is
$\text{gr}(\alpha_0)=\tilde{d}(S^3_p(9_{42}),[0])=d(S^3_1(9_{42}))=0$.
This proves the case of $m=1$.

Next, to compute $d(S^3_1(\overset{m}{\#}9_{42}))$, we consider $CFK^\infty(\overset{m}{\#}9_{42})\cong\otimes^mCFK^\infty(9_{42})$.
We denote the set of generators by $S_m=\{x_{i_1}\otimes x_{i_2}\otimes \cdots\otimes x_{i_m}\,|\,1\le i_k\le 9\}$ and $G_l^{(m)}=\mathbb{F}\langle U^{-l}\cdot x\,|\,x\in S_m\rangle$.
The complex $CFK^\infty(\overset{m}{\#}9_{42})$ is decomposed into the sum
$\oplus_{l\in\mathbb{Z}}G_l^{(m)}$ of the chain complexes.
Hence, we may consider each homology $H_\ast(G_l^{(m)}\cap \{\max(i,j)\ge 0\})$.

We define the tensor product $y\otimes \cdots\otimes y$ to be $y^{\otimes m}$.
Since the differential $\partial^\infty$ in $\otimes^mCFK^\infty(9_{42})$ is as follows:
$$\partial^\infty(z_1\otimes \cdots\otimes z_m)=\sum_{k=1}^mz_1\otimes\cdots \otimes\partial^\infty z_k\otimes \cdots\otimes z_m,$$
$H_\ast(G_l^{(m)})$ is generated by $U^{-l}\cdot\alpha^{\otimes m}$.
In fact, $\partial^\infty(\alpha^{\otimes m})=\sum_{k=1}^m\alpha\otimes\cdots\otimes \overset{k\text{-th}}{0}\otimes\cdots\otimes  \alpha=0$.
Since the generator $\alpha^{\otimes m}$ has the unique top grading in $G^{(m)}_l$, $\alpha^{\otimes m}$ does not lie in $\text{Im}(\partial^\infty)\cap G_l^{(m)}$.
Furthermore, since $HF^\infty(S^3_p(\overset{m}{\#}9_{42}))$ is standard, $H_\ast(G_i^{(m)})$ 
is isomorphic to $\mathbb{F}\langle U^{-l}\cdot\alpha^{\otimes m}\rangle\cong \mathbb{F}$.
This homology computes $H_\ast(G_{l}^{(m)}\{\max(i,j)\ge 0)\})$ for sufficiently large $l$.

We consider the generators of $H_\ast(G_0^{(m)}\{\max(i,j)\ge0\})$, i.e. the $l=0$ case.
Clearly, we have
$$G_0^{(m)}\{\max(i,j)\ge0\}=G_0^{(m)}\{\max(i,j)=0\}.$$
In this case, the chain complex $G_0^{(m)}\{\max(i,j)\ge0\}$ is isomorphic to 
$G_0^{(m)}\{i=0,j\le 0\}\oplus G_0^{(m)}\{j=0,i<0\}$.
The chain complex $G_0^{(m)}\{i=0,j\le 0\}$ is isomorphic to the $m$-times tensor product of the chain complex defined in Figure~\ref{Td},
i.e. this is isomorphic to the fundamental part of $CFK^\infty(\overset{m}{\#}3_1)$ with the top grading $0$.
Figure~\ref{G2} is the chain complex $G_0^{(2)}\{\max(i,j)\ge 0\}$ and
the sum of indicated classes in the picture presents the homological generator.
\begin{figure}[htbp]
\begin{center}
\unitlength 0.1in
\begin{picture}(  4.0800,  5.0500)( 13.9500,-10.0500)
%
\special{pn 20}%
\special{sh 1}%
\special{ar 1400 1000 10 10 0  6.28318530717959E+0000}%
%
\special{pn 20}%
\special{sh 1}%
\special{ar 1400 600 10 10 0  6.28318530717959E+0000}%
%
\special{pn 20}%
\special{sh 1}%
\special{ar 1800 1000 10 10 0  6.28318530717959E+0000}%
%
\special{pn 8}%
\special{pa 1740 1000}%
\special{pa 1460 1000}%
\special{fp}%
\special{sh 1}%
\special{pa 1460 1000}%
\special{pa 1528 1020}%
\special{pa 1514 1000}%
\special{pa 1528 980}%
\special{pa 1460 1000}%
\special{fp}%
%
\special{pn 8}%
\special{pa 1400 660}%
\special{pa 1400 940}%
\special{fp}%
\special{sh 1}%
\special{pa 1400 940}%
\special{pa 1420 874}%
\special{pa 1400 888}%
\special{pa 1380 874}%
\special{pa 1400 940}%
\special{fp}%
%
\special{pn 8}%
\special{pa 1400 500}%
\special{pa 1800 500}%
\special{ip}%
\end{picture}%
\caption{The chain complex of the fundamental part of $CFK^\infty(3_1)$.}
\label{Td}
\input{m=2.tex}
\caption{The chain complex $G_0^{(2)}\{\max(i,j)\ge 0\}$ and the homological generator.}
\label{G2}
\end{center}
\end{figure}

We denote the homological generators in $G_0^{(m)}\{i=0,j\le 0\}$ and $G_0^{(m)}\{j=0,i\le 0\}$ by
$\alpha_1$, and $\alpha_2$ respectively.
The classes $\alpha_i$ have $x_5^{\otimes m}$ as a non-zero component.

Here we claim the following:
\begin{lem}
\label{claim1}
$$H_\ast(G_0^{(m)}\{\max(i,j)\ge0\})\cong\mathbb{F}\cdot\beta,$$
where $\beta=(x_5+x_9)^{\otimes m}+(x_5+x_1)^{\otimes m}+x_5^{\otimes m}$.
The absolute grading of $\beta$ is $0$.
\end{lem}
\paragraph{Proof of Lemma \ref{claim1}.}
We show the element $\beta$ is a homological generator in the chain complex $G_0^{(m)}\{\max(i,j)\ge0\})$.

Since we have $\partial^\infty(x_5+x_9)=2x_6+x_4=x_4$, $\partial^\infty(x_5+x_1)=2x_4+x_6=x_6$, and $\partial^\infty x_5=x_6+x_4$, we get the following
$$\partial^\infty((x_5+x_9)^{\otimes m})=\sum_{k=1}^m(x_5+x_9)\otimes \cdots (x_5+x_9)\otimes \overset{k-\text{th}}{x_4}\otimes(x_5+x_9) \cdots\otimes (x_5+x_9)=\sum_{k=1}^mx_5\otimes \cdots\otimes x_5\otimes \overset{k-\text{th}}{x_4}\otimes x_5\cdots\otimes x_5$$
$$\partial^\infty((x_5+x_1)^{\otimes m})=\sum_{k=1}^m(x_5+x_1)\otimes \cdots (x_5+x_9)\otimes \overset{k-\text{th}}{x_6}\otimes(x_5+x_9) \cdots\otimes (x_5+x_9)=\sum_{k=1}^mx_5\otimes \cdots\otimes x_5\otimes \overset{k-\text{th}}{x_6}\otimes x_5\cdots\otimes x_5,$$
in $G_0^{(m)}\{\max(i,j)\ge 0\}$.
Summing these,
we get the following:
$$\partial^\infty((x_5+x_9)^{\otimes m}+(x_5+x_1)^{\otimes m}+x_5^{\otimes m})=0.$$
Thus, $\beta$ is a homological generator in $G_0^{(m)}\{\max(i,j)\ge 0\}$.

Any homological generator $z$ in $H_\ast(G_0^{(m)}\{\max(i,j)\ge0\})$ has the non-zero component $x_5^{\otimes m}$.
For, if $z$ does not have the non-zero component $x_5^{\otimes m}$,
then $z$ is presented by $y_1+y_2$, where $y_1\in G_0^{(m)}\{i=0,j< 0\}$ and $y_2\in G_0^{(m)}\{j=0,i< 0\}$.
However $y_1,y_2$ are not homological generators in $G_0^{(m)}\{i=0,j\le 0\}$ or $G_0^{(m)}\{j=0,i\le 0\}$ respectively.
Thus there exist $\bar{y}_1\in G_0^{(m)}\{i=0,j\le 0\}$ and $\bar{y}_2\in G_0^{(m)}\{j=0,i\le 0\}$
such that $y_1=\partial^{vert}\bar{y}_1$ and $y_2=\partial^{hor}\bar{y}_2$ respectively.
We may assume that $\bar{y}_i$ does not have $x_5^{\otimes m}$ as a non-zero component by reducing $\alpha_i$ if necessary.
This means $z=\partial^\infty(\bar{y}_1+\bar{y}_2)$ in $G_{0}^{(m)}\{\max(i,j)=0\}$.
This contradicts the fact that $z$ is a homological generator in $G_{0}^{(m)}\{\max(i,j)\ge0\}$.
Hence any homological generator in $G_0^{(m)}\{\max(i,j)\ge0\}$ has the non-zero component $x_5^{\otimes m}$. 
Since $z-\beta$ does not have $x_5^{\otimes m}$ as a non-zero component, it is not a homological generator.
This implies $z$ is homologous to $\beta$.

Therefore we have $H_\ast(G_0^{(m)}\{\max(i,j)\ge0\})\cong\mathbb{F}\cdot\beta$
and the grading of $\beta$ is $0$.
\hfill$\Box$\\

By Lemma~\ref{claim1}, $\beta$ is the minimal generator in $H_\ast(CFK^\infty(\overset{m}{\#}9_{42})\{\max(i,j)\ge 0\})$
coming from $HF^\infty(S^3_1(\overset{m}{\#}9_{42}))$.
Therefore we have $d(S^3_1(\overset{m}{\#}9_{42}))=\tilde{d}(S^3_{p}(\overset{m}{\#}9_{42}))=gr(\beta)=0$
\hspace{\fill} $\Box$ \\

\paragraph{Proof of Theorem \ref{thm3}.}

Since $\sigma(9_{42})= -2$ and the knot signature is additive, 
we have $\sigma(\overset{n+k}{\#}9_{42}) = -2(n+k)$.
Thus, by using Theorem \ref{thm2}, we have
\[
\gamma^0_{n\mathbb{C}P^2}(\overset{n+k}{\#}9_{42}) \geq \frac{-(-2(n+k))}{2} + 0 -n = k.
\]

We next construct a non-orientable surface 
$F_{n,k} \subset \text{punc}(n\mathbb{C}P^2)$ satisfying the following:

\begin{enumerate}
\item $\partial F_{n,k} = \overset{n+k}{\#}9_{42}$,
\item $\beta_1(F_{n,k})=k$, and 
\item $F_{n,k}$ represents zero in 
$
H_2 (\text{punc}(n\mathbb{C}P^2),
\partial(\text{punc}(n\mathbb{C}P^2)); \mathbb{Z}_2).
$
\end{enumerate}
\begin{figure}
\begin{center}
\includegraphics[clip]{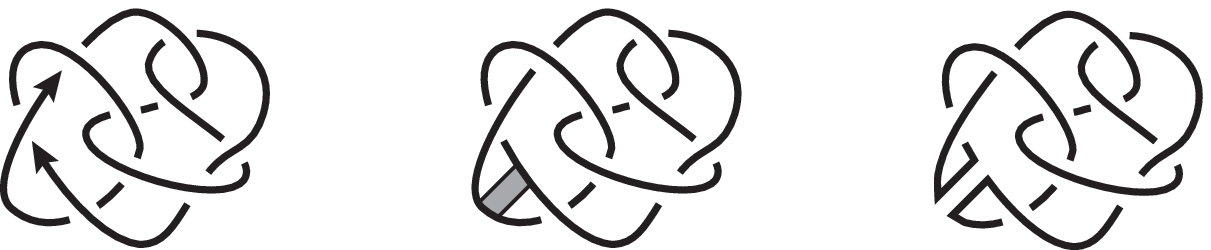}\\
\caption{$9_{42}$ bounds M\"{o}bius band in $B^4$}
\label{9_42 Mobius band}
\ \\
\includegraphics[clip]{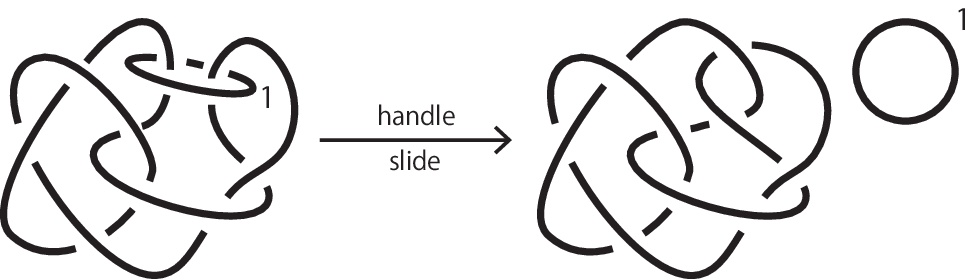}\\
\caption{$9_{42}$ bounds a disk in punc$\mathbb{C}P^2$}
\label{9_42 disk}
\end{center}
\end{figure}
The cobordisms in Figure \ref{9_42 Mobius band} and \ref{9_42 disk}
give a properly embedded M\"{o}bius band $M$ in $B^4$
with the boundary $9_{42}$, and a properly embedded disk $D$ in 
punc$\mathbb{C}P^2$ such that it bounds $9_{42}$
and represents zero in
$
H_2 (\text{punc}\mathbb{C}P^2,
\partial(\text{punc}\mathbb{C}P^2); \mathbb{Z}_2).
$
Taking the boundary connected sum of $n$ copies of  
$(\text{punc}\mathbb{C}P^2,D)$ and $k$ copies of $(B^4,M)$,
we have a new non-orientable surface $F_{n,k}$
satisfying the above properties from (1) to (3). This completes the proof.
\hspace{\fill} $\Box$ \\

\end{document}